\documentclass[leqno,intlimits,11pt]{amsart}

\usepackage[bookmarks=true]{hyperref}
\usepackage{amsmath}
\usepackage{amsthm,amssymb}
\usepackage[alphabetic]{amsrefs}
\usepackage{xypic}

\theoremstyle{plain} 

\theoremstyle{definition}
\newtheorem{thm}{Theorem}[section]
\newtheorem{lem}[thm]{Lemma}
\newtheorem{prop}[thm]{Proposition}

\newtheorem{cor}[thm]{Corollary}

\newtheorem{rmk}[thm]{Remark}

\newtheoremstyle{citing}
  {}
  {}
  {}
  {}
  {\bfseries}
  {.}
  {.5em }
  {\thmnote{#3}}

\theoremstyle{citing}
\newtheorem*{varthm}{} 

\numberwithin{equation}{section}

\newcommand{\Nguyen}{Nguy$\tilde{\hat{\text{e}}}$n}

\newcommand{\PP}{\mathbb{P}} 
\newcommand{\A}{\mathbb{A}} 
\newcommand{\CC}{\mathbb{C}} 
\newcommand{\ZZ}{\mathbb{Z}} 

\newcommand{\OO}{\mathcal{O}} 

\DeclareMathOperator{\Pic}{Pic} 
\DeclareMathOperator{\Hom}{Hom} 
\newcommand{\SU}{\mathcal{SU}} 
\newcommand{\U}{\mathcal{U}} 
\newcommand{\calC}{\mathcal{C}} 
\newcommand{\calD}{\mathcal{D}} 
\newcommand{\calI}{\mathcal{I}} 
\newcommand{\st}{\, : \,} 
\newcommand{\Thgen}{\Theta^{\text{gen}}} 
\newcommand{\bino}[2]{\bigl( \begin{smallmatrix} #1 \\ #2 \end{smallmatrix} \bigr)}

\begin{document}
\title[$\SU_X(3)$ and duality]{The moduli space of rank-3 vector bundles with trivial determinant
  over a curve of genus 2 and duality} \author{\Nguyen{} Quang Minh} \address{Department of
  Mathematics \\ University of Michigan \\ 525 East University Ave \\ Ann Arbor, MI 48109 \\ USA}
\email{ngminh@umich.edu} 

\begin{abstract}
  Let $\SU_X(3)$ be the moduli space of semi-stable vector bundles of rank 3 and trivial determinant
  on a curve $X$ of genus 2.  It maps onto $\PP^8$ and the map is a double cover branched over a
  sextic hypersurface called the Coble sextic.  In the dual $\PP^8$ there is a unique cubic
  hypersurface, the Coble cubic, singular exactly along the abelian surface of degree 1 line bundles
  on $X$.  We give a new proof that these two hypersurfaces are dual.  As an immediate corollary, we
  derive a Torelli-type result.
\end{abstract}

\maketitle

\section*{Introduction}

Let us fix once and for all a smooth projective curve $X$ of genus $g=2$.  We denote by $J^d(X)$, or
even by $J^d$ since we fixed $X$, the variety parametrizing classes of line bundles (or divisors) on
$X$ of degree $d$.  When $d=0$, we write $J$ for the Jacobian of $X$.  The variety $J^1$ carries a
canonical Riemann theta divisor
\begin{equation*}
\Theta = \{ L \in J^1 \st H^0(X,L) \neq 0 \}.
\end{equation*}
Moreover, we know that $3\Theta$ is very ample on $J^1$, so this give an embedding of $J^1$ into
$\PP^8=|3\Theta|^*$.  A. Coble in \cite{Cob17} shows that $J^1$ is set-theoretically cut out by 9
quadrics.  In \cite{Bar95}, W. Barth proves that this is even a scheme-theoretic intersection.  In
particular, if we denote by $\calI_{J^1}$ the ideal sheaf of $J^1$ in $\PP^8$, then
\begin{equation} \label{eq:9quadrics}
\dim H^0(\PP^8,\calI_{J^1}(2)) = 9,
\end{equation}
which can also be derived from the projective normality of the embedding of $J^1$ \cite{Koi76}.  It
turns out that the quadrics are the partial derivatives of a cubic, so there is a unique cubic
hypersurface $\calC_3$ singular exactly along the Abelian surface $J^1$.  The hypersurface $\calC_3$
has hence been dubbed the \emph{Coble cubic}.

Let $\SU_X(3)$ be the moduli space of semi-stable vector bundles of rank 3 and trivial determinant
on a curve $X$ of genus 2.  It maps onto $\PP^8=|3\Theta|$ and the map is a double cover branched
over a sextic hypersurface $\calC_6$.  This $\PP^8$ is the dual $\PP^8$ of the one in which
$\calC_3$ lies.  I. Dolgachev conjectured that $\calC_3$ and $\calC_6$ are dual varieties, and by
analogy with the case of the \emph{Coble quartic}, $\calC_6$ is known as the \emph{Coble sextic}.
Indeed, the Coble quartic, a quartic hypersurface in $\PP^7$, has an interpretation as the moduli
space $\SU_C(2)$ of semi-stable vector bundles of 2 and trivial determinant on a fixed
non-hyperelliptic curve $C$ of genus 3 (see \cite{NR87}).  Moreover, this quartic hypersurface is
singular exactly along the Kummer surface associated to $C$ (and embedded into $\PP^7$) and thanks
to its moduli space interpretation, C. Pauly \cite{Pau02} proved that the Coble quartic is
self-dual.

In this paper, we prove the following theorem:

\begin{varthm}[Theorem \ref{thm:duality}]
  The Coble hypersurfaces $\calC_3$ and $\calC_6$ are dual.
\end{varthm}

The result was first proved by A. Ortega Ortega \cite{Ort03} in her thesis.  We give here a
different proof, which uses a more thorough study of the variety $\calC_6$ and a more general
description of the dual map in terms of the vector bundles.  In particular we compute the degree of
its singular locus.  As a corollary, we derive a non-abelian Torelli result.

\section*{Acknowledgements}

I would like to express my sincere thanks to Igor Dolgachev for all the support, guidance and
encouragement during the research leading to this paper.  I would also like to acknowledge the
insightful discussions with Alessandro Verra, Angela Ortega Ortega, Mihnea Popa and Ravi Vakil.

\section{Definitions and preliminaries}

For a vector bundle $E$ or rank $n$ on $X$, we define its determinant
\begin{equation*}
\det(E) = \bigwedge^n E,
\end{equation*}
and its degree
\begin{equation*}
\deg(E) = \deg(\det(E)).
\end{equation*}
Then we define the slope of $E$ to be
\begin{equation*}
\mu(E) = \frac{\deg(E)}{n}.
\end{equation*}
We say that $E$ is stable (resp. semi-stable) if, for all proper subbundle $F$ of $E$, the
inequality
\begin{align*}
  \mu(F) &< \mu(E) & \text{( resp. } \mu(F) &\leq \mu(E) \text{ )}
\end{align*}
holds.

In \cite{NS64}, \cite{Ses67}, \cite{NR69}, moduli spaces of vector bundles are constructed.  For
some integers $n$ and $d$, we denote by $\U_X(n,d)$ the moduli space of semi-stable vector bundles
of rank $n$ and degree $d$ on $X$.  For a fixed line bundle $L \in J^d$ on $X$, we also denote by
$\SU_X(n,L)$ the moduli space of vector bundles of rank $n$ and determinant $L$.  For $L=\OO_X$, we
simply write $\SU_X(n)$.

We will focus our attention to $\SU_X(3)$ and recall a few facts proved in \cite{DN89}.  Much
analogously to the case of Jacobians, this moduli space carries a Cartier divisor represented by
\begin{equation*}
\Delta_{L} = \{ E \in \SU_X(3) \st H^0(X,E \otimes L) \neq 0 \}
\end{equation*}
for some fixed $L \in J^1$.  We will denote by
\begin{equation*}
\Thgen = \OO_{\SU_X(3)}(\Delta_L)
\end{equation*}
the corresponding line bundle.  $\Thgen$, which does not depend on the choice of $L$, is the ample
generator of $\Pic(\SU_X(3))$.  In \cite{BNR89}, the map $\Phi$ defined by the complete linear
system $|\Thgen|$ is given an alternate description: for $E \in \SU_X(3)$, let $D_E$ be
\begin{equation*}
D_E = \{ L \in J^1 \st H^0(X,E \otimes L) \neq 0 \}.
\end{equation*}
It is a divisor on $J^1$ from the linear system $|3\Theta|$.  This assignment defines an actual
morphism
\begin{equation*}
\xymatrix@R=0pt{D:\, \SU_X(3) \ar[r] & **[r] |3\Theta|, \\
  E \ar@{|->}[r] & **[r] D_E, }
\end{equation*}
which makes the following diagram commute
\begin{equation*}
\xymatrix@R=10pt{ & |\Thgen|^* \ar[dd]^{\cong} \\
  \SU_X(3) \ar[ru]^{\Phi} \ar[rd]_{D} &  \\
  & |3\Theta| \ .}
\end{equation*}
Therefore in the following we will identity $\Phi$ and $D$.  It is known that $\Phi$ is a finite map
of degree 2.  A first unpublished proof was given by D. Butler and I. Dolgachev using the Verlinde
formula, but another beautiful proof can be found in \cite{Las96}.  The branch divisor (which is
isomorphic to the ramification locus) is a hypersurface of degree 6 in $|3\Theta| \cong \PP^8$.  We
denote it by $\calC_6$.  Since $\PP^8$ is smooth, the singular locus
\begin{equation*}
\Sigma = \text{Sing}(\calC_6) \subset \PP^8
\end{equation*}
is exactly the singular locus $\Sigma'$ of $\SU_X(3)$ corresponding to strictly semi-stable vector
bundles (i.e. semi-stable but not stable.)  We will keep the two notations, $\Sigma$ and $\Sigma'$,
in order to make clear in what space we are.  Let us recall that $\Sigma$ is of dimension 5.  There
is a normalization map
\begin{equation*}
\xymatrix@R=0pt{\nu:\, \U_X(2,0) \ar[r] & **[r] \Sigma', \\
  F \ar@{|->}[r] & **[r] F \oplus \det(F)^*. }
\end{equation*}
It is a proper birational map.

Finally, for any variety $Z$ and coherent sheaf $\mathcal{F}$ on $Z$, we will write
\begin{equation*}
h^i(Z,\mathcal{F}) = \dim H^i(Z,\mathcal{F}).
\end{equation*}

\section{The degree of $\Sigma$}

By definition, the degree of $\Sigma$ is
\begin{equation*}
\deg(\Sigma) = \int_{\PP^8} \Sigma \cdot (H^5),
\end{equation*}
where $H$ is the class of a hyperplane in $\PP^8$.

\begin{prop} \label{prop:degSigma1}
  The degree of $\Sigma$ is given by
\begin{equation*}
\deg(\Sigma) = \int_{U_X(2,0)} (\nu^*(\Thgen)^5).
\end{equation*}
\end{prop}

\begin{proof}
  Since $\Sigma = \Phi_* \Sigma'$, we apply the projection formula and see that
\begin{equation*}
(\ \Sigma \cdot (H^5)\ ) = (\ \Phi_* \Sigma' \cdot (H^5)\ ) = (\ \Sigma' \cdot (\Phi^* H)^5\ ) = (\ 
\Sigma' \cdot (\Thgen)^5\ ).
\end{equation*}
To compute the last intersection product, we use the normalization map $\nu$:
\begin{equation*}
\int_{\SU_X(3)} \Sigma' \cdot (\Thgen)^5 = \int_{\SU_X(3)} (\ \nu_* U_X(2,0) \cdot (\Thgen)^5 \ ) =
\int_{U_X(2,0)} (\nu^*(\Thgen)^5).
\end{equation*}
\end{proof}

Recall that a divisor representing $\Thgen$ is
\begin{equation*}
\Delta_{L} = \{ E \in \SU_X(3) \st H^0(X,E \otimes L) \neq 0 \}
\end{equation*}
for some $L \in J^1$.  Then
\begin{equation} \label{eq:support}
\begin{split}
  \nu^*(\Delta_{L}) &= \{ F \in \U_X(2,0) \st H^0(X,(F \otimes L) \oplus (\det(F)^* \otimes L)) \neq 0 \}, \\
  &= \{ F \in \U_X(2,0) \st H^0(X,F\ \otimes L) \neq 0 \} \\
  & \hspace{2cm} \cup \{ F \in \U_X(2,0) \st H^0(X,\det(F)^* \otimes L) \neq 0) \}.
\end{split}
\end{equation}
To deal with this, let us recall that the determinant map
\begin{equation*}
\xymatrix@R=0pt{\det:\, \U_X(2,0) \ar[r] & **[r] J, \\
  F \ar@{|->}[r] & **[r] \det(F), }
\end{equation*}
is actually a $\PP^3$-fibration.  Indeed, for $a \in J$, the fiber over $a$ is $\SU_X(2,a) \cong
\PP^3$.  So for a fixed $L \in J^1$, we define the map $\pi$ as the following composition:
\begin{equation*}
\xymatrix{ \U_X(2,0) \ar[r]^{\pi} \ar[d]_{\det} & J^1 \\
  J \ar[r]^{-1}_{\cong} & J \ar[u]_{\otimes L}^{\cong} }
\end{equation*}
and we see that $\pi$ is also a $\PP^3$-fibration.  Let us also recall that $\U_X(2,0)$ has a
generalized theta divisor
\begin{equation*}
\Thgen_{\U} = \OO_{\U_X(2,0)}(\Delta'_L)
\end{equation*}
where
\begin{equation*}
\Delta'_L = \{ F \in U_X(2,0) \st H^0(X,F \otimes L) \neq 0 \} \subset \U_X(2,0),
\end{equation*}
still for some $L \in J^1$.  So we see that at the divisorial level (or set-theoretically)
\begin{equation*}
\nu^*(\Delta_L) = \Delta'_L \cup \pi^*(\Theta),
\end{equation*}
and as an isomorphism class of line bundles on $\U_X(2,0)$,
\begin{equation} \label{eq:nuThgen}
\nu^*(\Thgen) = \Thgen_{\U} + \pi^*(\Theta).
\end{equation}

\begin{prop} \label{prop:degSigma}
  The degree of the singular locus $\Sigma$ of $\calC_6$ in $\PP^8$ is
\begin{equation*}
\deg(\Sigma) = 45.
\end{equation*}
\end{prop}

Here is the motivation of the proof, which will follow from three lemmas.  Putting
\eqref{eq:nuThgen} and Proposition \ref{prop:degSigma1} together, we obtain
\begin{align*}
  \deg(\Sigma) &= \int_{\U_X(2,0)} [\nu^*(\Thgen)]^5 = \int_{\U_X(2,0)} \left([\Thgen_{\U}] + [\pi^*(\Theta)]\right)^5 \\
  &= \int_{\U_X(2,0)} [\Thgen_{\U}]^5 \ + \ \bino{5}{1} \int_{\U_X(2,0)} [\Thgen_{\U}]^4 \cdot [\pi^*(\Theta)] \\
  & \quad + \bino{5}{2} \int_{\U_X(2,0)} [\Thgen_{\U}]^3 \cdot [\pi^*(\Theta)]^2 \ + \ \bino{5}{3} \int_{\U_X(2,0)} [\Thgen_{\U}]^2 \cdot [\pi^*(\Theta)]^3 \\
  & \quad + \bino{5}{4} \int_{\U_X(2,0)} [\Thgen_{\U}] \cdot [\pi^*(\Theta)]^4 \ + \ 
  \int_{\U_X(2,0)} [\pi^*(\Theta)]^5.
\end{align*}
But $\pi$ is a fibration and by the projection formula, some intersection cycles are zero:
\begin{equation*}
\pi^*(\Theta) \cdot \pi^*(\Theta) \cdot \pi^*(\Theta) = 0,
\end{equation*}
hence
\begin{multline} \label{eq:degSigma2}
  \deg(\Sigma) = \int_{\U_X(2,0)} [\Thgen_{\U}]^5 \ + \ 5 \int_{\U_X(2,0)} [\Thgen_{\U}]^4 \cdot [\pi^*(\Theta)] \\
  + 10 \int_{\U_X(2,0)} [\Thgen_{\U}]^3 \cdot [\pi^*(\Theta)]^2.
\end{multline}
There are essentially three terms in the sum and we will treat each one separately.

\begin{lem} 
\begin{equation*}
\int_{\U_X(2,0)} [\Thgen_{\U}]^5 = 5.
\end{equation*}
\end{lem}

\begin{proof}
  For this we use the following \'etale covering:
\begin{equation} \label{eq:psi}
\xymatrix@R=0pt{**[l] \psi:\, \SU_X(2) \times J  \ar[r]^{16:1} & **[r] \U_X(2,0), \\
**[l] (F,L) \ar@{|->}[r] & **[r] F \otimes L. }
\end{equation}
Notice in this case that $\det(F\otimes L) = L^2$.  We know that
\begin{equation} \label{eq:pullbackThgen}
\psi^*(\Thgen_{\U}) = pr_1^*(\OO_{\PP^3}(1)) \otimes pr_2^*(2\Theta_0),
\end{equation}
where $pr_1$ (resp. $pr_2$) is the projection from $\SU_X(2) \times J$ onto the first (resp. second)
factor, and $\Theta_0$ is a (symmetric) theta divisor on $J$, a principal polarization and
\begin{equation*}
\int_{J} \Theta_0 \cdot \Theta_0 = 2.
\end{equation*}
So when we apply the projection formula, we get
\begin{align*}
  [\psi^*(\Thgen_{\U})]^5 &= ([pr_1^*(\OO_{\PP^3}(1))] + [pr_2^*(2\Theta_0)])^5, \\
  16 \int_{\U_X(2,0)} [\Thgen_{\U}]^5 &= \bino{5}{2} \quad \int_{\SU_X(2) \times J} [pr_1^*(\OO_{\PP^3}(1))]^3 \cdot [pr_2^*(2\Theta_0)]^2, \\
  16 \int_{\U_X(2,0)} [\Thgen_{\U}]^5 &= 10 \times 1 \times 8 = 80.
\end{align*}
Then, by dividing by 16, we get the result.
\end{proof}

\begin{lem}
\begin{equation*}
\int_{\U_X(2,0)} [\Thgen_{\U}]^4 \cdot [\pi^*(\Theta)] = 4.
\end{equation*}
\end{lem}

\begin{proof}
  We will again pull back the intersection by the \'etale map $\psi$ \eqref{eq:psi}.  Recall
  \eqref{eq:support} that the class $\pi^*(\Theta)$ can be represented by the divisor
\begin{equation*}
D = \{ F \in \U_X(2,0) \st H^0(X,\det(F)^* \otimes L) \neq 0 \},
\end{equation*}
where $L$ is our fixed line bundle on $X$ of degree 1.  Then, pulling back $D$:
\begin{align*}
  \psi^*(D) &= \{ (F,N) \in \SU_X(2) \times J \st H^0(X,\det(F)^* \otimes N^{-2} \otimes L) \neq 0 \}, \\
  &= \{ (F,N) \in \SU_X(2) \times J \st H^0(X,N^{-2} \otimes L) \neq 0 \},
\end{align*}
since $\det(F)=\OO_X$.  Therefore its class is equal to
\begin{equation*}
[\psi^*(\pi^*(\Theta))] = pr_2^* \ [-2]^* \ t_L^*([\Theta]),
\end{equation*}
where $[-2]$ is the multiplication by 2 on $J$ and $t_L:\ J \to J^1$ is the translation induced by
tensoring by $L$.  This can be easily seen from the following commutative diagram:
\begin{equation*}
\xymatrix{ **[l] \SU_X(2) \times J \ar[r]^{\psi} \ar[d]_{pr_2} & \U_X(2,0) \ar[r]^{\pi} \ar[d]^{\det} & J^1 \\
  J \ar[r]^{[2]} & J \ar[r]^{[-1]} & **[r] J \quad . \ar[u]_{t_L=\ \otimes L} }
\end{equation*}
First, $t_L^*(\Theta)$ is just a translate of a symmetric theta divisor $\Theta_0$ of $J$, so for
intersection purposes, we might as well assume that it is $\Theta_0$.  Therefore
\begin{equation*}
[-2]^*([\Theta_0]) = [2]^*([\Theta_0]) = [k\Theta_0],
\end{equation*}
for some number $k$ to be determined.  Since $[2]$ is a endomorphism of $J$ of degree 16, we have on
the one hand,
\begin{equation*} 
\int_{J} ([2]^*([\Theta_0]))^2 = 16 \int_J [\Theta_0]^2 = 32.
\end{equation*}
On the other hand,
\begin{equation*}
\int_{J} ([2]^*([\Theta_0]))^2 = \int_J [k\Theta_0]^2 = 2k^2,
\end{equation*}
so $k=4$.  Thus
\begin{equation*}
[\psi^*(\pi^*(\Theta))] = pr_2^*([4\Theta_0]).
\end{equation*}
Finally, using \eqref{eq:pullbackThgen} again:
\begin{align*}
  \psi^*[\Thgen_{\U}]^4 \cdot \psi^*[\pi^*(\Theta)] &= \left([pr_1^*(\OO_{\PP^3}(1))] + [pr_2^*(2\Theta_0)]\right)^4 \cdot [pr_2^*(4\Theta_0)], \\
  16 [\Thgen_{\U}]^4 \cdot [\pi^*(\Theta)] &= \bino{4}{3} [pr_1^*(\OO_{\PP^3}(1))]^3 \cdot
  [pr_2^*(2\Theta_0)] \cdot [pr_2^*(4\Theta_0)].
\end{align*}
Hence
\begin{equation*}
\int_{\U_X(2,0)} [\Thgen_{\U}]^4 \cdot [\pi^*(\Theta)] = \frac{1}{16} \times 4 \times 1 \times 16 =
4.
\end{equation*}
\end{proof}

\begin{lem}
\begin{equation*}
\int_{\U_X(2,0)} [\Thgen_{\U}]^3 \cdot [\pi^*(\Theta)]^2 = 2.
\end{equation*}
\end{lem}

\begin{proof}
  This can be calculated directly and is fairly easy.  First, notice that
\begin{equation*}
[\pi^*(\Theta)]^2 = \pi^*([\Theta]^2) = \pi^*([\text{2 points}]) = [\text{2 fibers}].
\end{equation*}
We know however that the fibers of $\pi$ are of the form $\SU_X(2,a)$ but we know \cite{NR69} that
$\SU_X(2,a) \cong \PP^3$, for $a \in J$, so
\begin{equation*}
[\Thgen_{\U}] \cdot [\text{fiber}] = [\Thgen_{\SU_X(2,a)}] = [\OO_{\PP^3}(1)]
\end{equation*}
as a cycle of $\SU_X(2,a) \cong \PP^3$.  Therefore
\begin{equation*}
\int_{\U_X(2,0)} [\Thgen_{\U}]^3 \cdot [\text{2 fibers}] = 2 \int_{\PP^3} [\OO_{\PP^3}(1)]^3 = 2.
\end{equation*}
\end{proof}

\begin{proof}[Proof of Proposition \ref{prop:degSigma}]
  It suffices now to collect the numbers obtained in the previous three lemmas and put them into the
  equation \eqref{eq:degSigma2}:
\begin{align*}
  \deg(\Sigma) &= \int_{\U_X(2,0)} [\Thgen_{\U}]^5 \ + \ 5 \int_{\U_X(2,0)} [\Thgen_{\U}]^4 \cdot [\pi^*(\Theta)] \\
  & \hspace{4cm} + 10 \int_{\U_X(2,0)} [\Thgen_{\U}]^3 \cdot [\pi^*(\Theta)]^2, \\
  &= 5 + 5 \times 4 + 10 \times 2 = 45.
\end{align*}
\end{proof}

\section{A map given by quadrics}

Let $P_3$ be the homogeneous cubic polynomial defining the Coble cubic $\calC_3$.  The motivation
here is to interpret the dual map, given by quadrics,
\begin{equation*}
\xymatrix@R=0pt{
  \calD:\ \calC_3 \subset |3\Theta|^* \ar@{-->}[r] & **[r] \check{\calC}_3 \subset |3\Theta| \\
  p \ar@{|->}[r] & **[r] T_p(\calC_3) = \left[ \frac{\partial P_3}{\partial X_i}
  \right]_{i=0,\dots8} }
\end{equation*}
in terms of vector bundles.  So we are trying to construct (i.e. identify) a rational map
\begin{equation*}
\Psi:\ \calC_3 \dasharrow \SU_X(3).
\end{equation*}

Recall that
\begin{equation*}
X \cong \Theta = \{ L \in J^1 \st H^0(X,L) \neq 0 \}.
\end{equation*}
Now for every $a \in J$, we write
\begin{equation*}
X_a = \Theta + a \subset J^1.
\end{equation*}
Recall also that $J^1$ is embedded by complete linear series into $|3\Theta|^*$.  Let us restrict
the embedding to $X_a$: we write $\vartheta_a = \OO_{X_a}(\Theta|_{X_a})$, a line bundle of degree 2
on $X_a$.  It is important to notice that if $a \neq b$ in $J$, then $\vartheta_a \ncong
\vartheta_b$.  By Riemann-Roch, we get
\begin{equation*}
h^0(X_a,\vartheta_a{}^3) = 5
\end{equation*}
and we denote by $\PP^4_a = |\vartheta_a{}^3|^*$ the linear span of $X_a$ in $\PP^8 = |3\Theta|^*$.
The goal will be to define a rational map on $\PP^4_a$ and before we go any further, we need the
following result.

\begin{prop}
  The span of $X_a$ in $|3\Theta|^*$ lies in $\calC_3$:
\begin{equation*}
\PP^4_a \subset \calC_3.
\end{equation*}
\end{prop}

\begin{proof}
  Suppose $\PP^4_a \nsubseteq \calC_3$, then $\PP^4_a \cap \calC_3 = V_3$ a cubic threefold of
  $\PP^4_a$.  Since $\calC_3$ is singular exactly along $J^1$, then $X_a \subset \text{Sing}(V_3)$.
  Let $\ell$ be a secant line to $X_a$, so
\begin{equation*}
\int_{\PP^4_a} \ell \cdot V_3 \geq 4
\end{equation*}
which implies that $\ell \subset V_3$.  Therefore the secant variety $\text{Sec}(X_a)$, a threefold
since $X_a$ does not lie in a plane, is contained in $V_3$.  But we will show that
\begin{equation*}
\deg \text{Sec}(X_a) = 8,
\end{equation*}
a contradiction.  Indeed, let $\ell$ be a general line in $\PP^4_a$, it intersects $\text{Sec}(X_a)$
at $d$ points.  Therefore, when we project from $\ell$, $X_a$ is mapped to a plane sextic curve of
geometric genus 2 with $d$ nodes.  Since the arithmetic genus of a plane sextic curve is 10, we see
that $d=8$.
\end{proof}

Let $x \in \PP^4_a - X_a$, it corresponds to a hyperplane $V_x$ in $H^0(X_a,\vartheta_a{}^3)$, i.e
it is equivalent to giving a linear form on $H^0(X_a,\vartheta_a{}^3)$:
\begin{equation*}
0 \xrightarrow{} V_x \xrightarrow{j_x} H^0(X_a,\vartheta_a{}^3) \xrightarrow{x} \CC \xrightarrow{}
0.
\end{equation*}
Since $\vartheta_a{}^3$ is very ample on $X_a$, it is generated by its global sections.  But $x
\notin X_a$, so $V_x$ still generates $\vartheta_a{}^3$.  We write down the evaluation (exact)
sequences:
\begin{equation} \label{eq:commdiag}
\xymatrix{
0 \ar[r] & E_x \ar[r] \ar[d]_i & V_x \otimes \OO_{X_a} \ar[r]^{e_x} \ar[d]_{j_x} & \vartheta_a{}^3 \ar[r] \ar@{=}[d] & 0 \\
0 \ar[r] & M \ar[r] \ar[d] & H^0(X_a,\vartheta_a{}^3) \otimes \OO_{X_a} \ar[r]^{\qquad e}\ar[d]_{x} & \vartheta_a{}^3 \ar[r] & 0 \\
& \OO_{X_a} \ar@{=}[r] & \OO_{X_a} & & }
\end{equation}
where $i$ is an inclusion and the lower row comes from the snake lemma.  The sheaves $E_x$ and $M$
are locally free so we will see them as vector bundles of rank 3 and 4 respectively.  We then see
that they both have degree $-6$.  So, their slopes are
\begin{align*}
  \mu(E_x) &= -2, & \mu(M) &= -\frac{3}{2}.
\end{align*}

\begin{lem}
  The vector bundle $E_x$ is semi-stable.
\end{lem}

\begin{proof}
  Suppose $F$ is a subbundle of $E_x$.  It is also a subbundle of $M$, but $M$ is known to be stable
  \cite{EL92} because $\deg(\vartheta_a{}^3)=6$.  So $\mu(F) < \mu(M) = -\frac{3}{2}$, therefore
  $\mu(F) \leq -2$ because $F$ is of rank 1 or 2, i.e $\mu(F) \leq \mu(E_x)$.
\end{proof}

In particular, $E_x(\vartheta_a)$ is semi-stable (because $E_x$ is) of rank 3 and fits in the
twisted evaluation sequence
\begin{equation} \label{eq:ses}
0 \to E_x(\vartheta_a) \to V_x \otimes \vartheta_a \to \vartheta_a{}^4 \to 0.
\end{equation}
Therefore,
\begin{equation*}
\det(E_x(\vartheta_a)) = \OO_{X_a}.
\end{equation*}
So this assignment defines a rational map $\Psi$ from $\PP^4_a$ to $\SU_X(3)$, regular outside of
$X_a$:
\begin{equation*}
\xymatrix@R=0pt{
  \Psi:\ \PP^4_a-X_a \ar[r] & **[r] \SU_X(3) \\
  \qquad x \ar@{|->}[r] & **[r] \Psi(x) = E_x(\vartheta_a). }
\end{equation*}

We will now study this map to see that it is defined by quadrics.  Let us twist the exact sequence
\eqref{eq:ses} by $\omega_{X_a} \otimes \vartheta_a^{-1}$, we get:
\begin{equation*}
0 \to E_x(\omega_{X_a}) \to V_x \otimes \omega_{X_a} \to \vartheta_a{}^3 \otimes \omega_{X_a} \to 0.
\end{equation*}
Then, by applying Riemann-Roch, we find that
\begin{align*}
  h^0(X_a,V_x \otimes \omega_{X_a}) &= 8, & h^0(X_a,\vartheta_a{}^3 \otimes \omega_{X_a}) &=7,
\end{align*}
so
\begin{equation*}
h^0(X_a,E_x(\omega_{X_a})) \geq 1,
\end{equation*}
which means that $E_x(\omega_{X_a})$, of degree 0, has sections, i.e. there is a non trivial
morphism
\begin{equation*}
\OO_{X_a} \to E_x(\omega_{X_a}).
\end{equation*}
Because, the two vector bundles are of degree 0, this morphism has to be an injection and the
quotient is also a vector bundle.  Now twisting back by $\vartheta_a \otimes \omega_{X_a}^{-1}$, we
obtain the following short exact sequence of vector bundles:
\begin{equation*}
0 \to \vartheta_a \otimes \omega_{X_a}^{-1} \to E_x(\vartheta_a) \to G \to 0.
\end{equation*}
Because these vector bundles are all of degree 0, it follows that
\begin{equation*}
E_x(\vartheta_a) \ \sim_S \ (\vartheta_a \otimes \omega_{X_a}^{-1}) \oplus G,
\end{equation*}
i.e. they represent the same class in the moduli space $\SU_X(3)$.  Moreover, we see that
\begin{align*}
  \det G &= \omega_{X_a} \otimes \vartheta_a^{-1}, & \nu(G) = E_x(\vartheta_a).
\end{align*}
Recall that $\SU_X(2,\omega_{X_a} \otimes \vartheta_a^{-1})$ sits naturally in $\U_X(2,0)$ as a
fiber of the determinant map.  It is also easy to check that
\begin{equation*}
\nu|_{\SU_X(2,\omega_{X_a} \otimes \vartheta_a^{-1})}:\ \SU_X(2,\omega_{X_a} \otimes
\vartheta_a^{-1}) \to \Sigma' \subset \SU_X(3)
\end{equation*}
is isomorphic onto its image and that the composition
\begin{equation} \label{eq:Dnu}
D \circ \nu:\ \SU_X(2,\omega_{X_a} \otimes \vartheta_a^{-1}) \to \Sigma' \to \Sigma \subset |3\Theta|=\PP^8
\end{equation}
embeds $\SU_X(2,\omega_{X_a} \otimes \vartheta_a^{-1})$ as a linear space, i.e $\PP^3$, in
$|3\Theta|$.  So we will write $\PP^3_a$ for $\SU_X(2,\omega_{X_a} \otimes \vartheta_a^{-1})$ and
see it as a subspace of $\Sigma$ or $\Sigma'$ interchangeably.

Therefore we just proved that the map $\Psi$ actually lands into $\PP^3_a$.

\begin{prop} \label{prop:Psi}
  The rational map
\begin{equation*}
\xymatrix@R=0pt{
  \Psi:\ \PP^4_a \ar@{-->}[r] & **[r] \PP^3_a \\
  x \ar@{|->}[r] & **[r] \Psi(x) = E_x(\vartheta_a) }
\end{equation*}
is given by a linear system of quadrics.
\end{prop}

\begin{proof}
  The degree of the linear system defining $\Psi$ is the degree of $\Psi^*(\Thgen)$.  At the
  divisorial level, if we fix $L \in J^1$, this is just
\begin{equation*}
\Psi^*(\Delta_L) = \{ x \in \PP^4_a \st H^0(X_a,E_x(\vartheta_a)\otimes L) \neq 0 \}
\end{equation*}
where $E_x(\vartheta_a) = \Psi(x)$ (and recall that $\Thgen=\OO_{\SU_X(3)}(\Delta_L)$).  Let us
choose $L \in J^1$ so that $\vartheta_a \otimes L$ is globally generated.  Such an $L$ exists
because, by Riemann-Roch, global generation is equivalent to
$h^0(X_a,\omega_{X_a}\otimes\vartheta_a^{-1}\otimes L^{-1}(p))=0$ for all $p\in X_a$.  The
assignment
\begin{equation*}
p \to \omega_{X_a} \otimes \vartheta_a^{-1} \otimes L^{-1}(p)
\end{equation*}
embeds $X_a$ into $J$, and the choice of $L$ is just the choice of a translation.  Certainly, we can
choose $L$ so that the image of $X_a$ in $J$ avoids the origin.  Next, if we twist the commutative
exact diagram \eqref{eq:commdiag} by $\vartheta_a \otimes L$, we get the following commutative "long
exact" diagram:
\begin{equation} \label{eq:commdiag2}
\xymatrix{& 0 \ar[d] & 0 \ar[d] & \\
0 \ar[r] & H^0(E_x(\vartheta_a)\otimes L) \ar[r] \ar[d] & V_x \otimes H^0(\vartheta_a \otimes L) \ar[r]^{\quad e_x} \ar[d] & H^0(\vartheta_a{}^4 \otimes L) \ar@{=}[d] \\
0 \ar[r] & H^0(M(\vartheta_a)\otimes L) \ar[r] \ar[d] & H^0(\vartheta_a{}^3) \otimes H^0(\vartheta_a \otimes L) \ar[d]_{g(x)} \ar[r]^{\qquad \quad e} & H^0(\vartheta_a{}^4 \otimes L) \\
& H^0(\vartheta_a \otimes L) \ar@{=}[r] & H^0(\vartheta_a{} \otimes L)  & }
\end{equation}
where the cohomology groups are taken over $X_a$ and where the map $g(x)$ can be described as
follows.  The contraction
\begin{equation*}
H^0(\vartheta_a{}^3) \otimes H^0(\vartheta_a{}^3)^* \otimes H^0(\vartheta_a \otimes L) \to
H^0(\vartheta_a \otimes L)
\end{equation*}
defines a linear map
\begin{equation*}
g:\ H^0(\vartheta_a{}^3)^* \to \Hom (H^0(\vartheta_a{}^3) \otimes H^0(\vartheta_a \otimes L),
H^0(\vartheta_a \otimes L)),
\end{equation*}
and we will abuse notation and often write $x$ for both an element of $\PP H^0(\vartheta_a{}^3)^*$
and any of its representatives in $H^0(\vartheta_a{}^3)^*$.  So for all $x \in
H^0(\vartheta_a{}^3)^*$, recall that we write $V_x = \ker(x)$, and
\begin{gather*}
  g(x):\ H^0(\vartheta_a{}^3) \otimes H^0(\vartheta_a \otimes L) \to H^0(\vartheta_a \otimes L), \\
  \ker(g(x)) = V_x \otimes H^0(\vartheta_a \otimes L).
\end{gather*}

Now, since $\dim(V_x \otimes H^0(\vartheta_a \otimes L)) = \dim(H^0(\vartheta_a{}^4 \otimes L)) =
8$, we have exactly
\begin{equation*}
\Psi^*(\Delta_L) = \{ x \in \PP^4_a \st \text{the map $e_x$ degenerates} \}.
\end{equation*}
Since $h^0(\vartheta_a \otimes L)=2$ and $\vartheta_a \otimes L$ is globally generated, we apply the
base point free pencil trick (see for instance \cite{ACGH85}) and see that the map $e$ has a
two-dimensional kernel
\begin{equation*}
\ker(e) = H^0(\vartheta_a{}^2 \otimes L^*),
\end{equation*}
so by restricting $g$, we get a map
\begin{equation*}
g':\ H^0(\vartheta_a{}^3)^* \to \Hom(H^0(\vartheta_a{}^2 \otimes L^*),H^0(\vartheta_a \otimes L))
\cong \Hom(\CC^2,\CC^2),
\end{equation*}
and we can rewrite the commutative exact diagram \eqref{eq:commdiag2} as follows:
\begin{equation} \label{eq:commdiag3}
\xymatrix{& 0 \ar[d] & 0 \ar[d] & \\
0 \ar[r] & H^0(E_x(\vartheta_a)\otimes L) \ar[r] \ar[d] & V_x \otimes H^0(\vartheta_a \otimes L) \ar[r]^{e_x} \ar[d] & H^0(\vartheta_a{}^4 \otimes L) \ar@{=}[d] \\
0 \ar[r] & H^0(\vartheta_a{}^2 \otimes L^*) \ar[r] \ar[d]_{g'(x)} & H^0(X_a,\vartheta_a{}^3) \otimes H^0(\vartheta_a \otimes L) \ar[d]_{g(x)} \ar[r]^{\qquad \quad e} & H^0(\vartheta_a{}^4 \otimes L) \\
& H^0(\vartheta_a{} \otimes L) \ar@{=}[r] & H^0(\vartheta_a{} \otimes L). & }
\end{equation}
It is then clear that $e_x$ degenerates, i.e. $H^0(E_x(\vartheta_a)\otimes L) \neq 0$, exactly when
$g'(x)$ degenerates, since $h^0(\vartheta_a{}^2 \otimes L^*) = h^0(\vartheta_a \otimes L) = 2$.  So
\begin{equation} \label{eq:g'(x) degenerates}
\Psi^*(\Delta_L) = \{ x \in \PP H^0(\vartheta_a{}^3)^* \st \det(g'(x)) = 0 \},
\end{equation}
i.e. it is the pull-back under $g'$ of the discriminant locus of
\begin{equation*}
\Hom(H^0(\vartheta_a{}^2 \otimes L^*),H^0(\vartheta_a \otimes L)),
\end{equation*}
a quadric.
\end{proof}

\section{Restriction of the dual map}

Now that we have proved that the map $\Psi$ is given by quadrics, we would like to show that it is
the restriction of the dual map $\calD:\ \calC_3 \dasharrow \check{\calC}_3$.  But first, let us
study the restriction of the dual map.  Recall that $X \cong X_a = \Theta+a \subset J^1$, for $a \in
J$.  Moreover, $\PP^4_a$ is the linear span of $X_a$ in $|3\Theta|^*=\PP^8$.  The dual map $\calD$
is given by the linear system of quadrics in $\PP^8$ containing $J^1$, so when restricted to
$\PP^4_a$, it is given by quadrics in $\PP^4_a$ containing $X_a$.

\begin{prop}
  The restriction of the dual map
\begin{equation*}
\calD|_{\PP^4_a} :\ \PP^4_a \dasharrow |3\Theta|
\end{equation*}
is given by the complete linear system $|\calI_{X_a}(2)|$ of quadrics in $\PP^4_a$ containing $X_a$.
\end{prop}

\begin{proof}
  The statement of the proposition is equivalent to the following:
\begin{equation*}
|\calI_{X_a}(2)|^* \subset |\calI_{J^1}(2)|^*
\end{equation*}
where $\calI_{J^1}$ denotes the sheaf of ideals of $J^1$ in $\PP^8$, or equivalently the natural
restriction map
\begin{equation} \label{eq:alpha}
\alpha:\ |\calI_{J^1}(2)| \to |\calI_{X_a}(2)|
\end{equation}
is surjective.  We have the exact sequence
\begin{equation*}
0 \to \calI_{X_a} \to \OO_{\PP^4_a} \to \OO_{X_a} \to 0,
\end{equation*}
which, twisted by $\OO_{\PP^4_a}(2)$, gives
\begin{equation} \label{eq:seqXa}
0 \to \calI_{X_a}(2) \to \OO_{\PP^4_a}(2) \to \OO_{X_a}(\vartheta_a{}^6) \to 0.
\end{equation}
It is known that $\calI_{X_a}$ is 3-regular because $X_a$ is neither rational nor elliptic (see for
instance \cite{GLP83}), so $H^1(\PP^4_a, \calI_{X_a}(2))=0$, and from the long exact sequence
associated to \eqref{eq:seqXa} we get
\begin{equation*}
\dim H^0(\PP^4_a,\calI_{X_a}(2)) = 4.
\end{equation*}
Let us recall \eqref{eq:9quadrics} that
\begin{equation*}
\dim H^0(\PP^8,\calI_{J^1}(2)) = 9.
\end{equation*}
Our goal is to prove that the map $\alpha$ \eqref{eq:alpha} is surjective.  Assume that it is not,
i.e. that the rank of $\alpha$ is strictly less than 4.  Then we can choose a basis of
$H^0(\PP^8,\calI_{J^1}(2))$ that consists of at least 6 quadrics that contain $\PP^4_a$ and at most
3 quadrics that do not.  Let us call $\mathcal{B}$ the base locus restricted to $\PP^4_a$ of the
latter quadrics:
\begin{equation*}
\mathcal{B} \subset \PP^4_a.
\end{equation*}
Then $\mathcal{B}$ is of dimension at least 1.  We see that the base locus of
$H^0(\PP^8,\calI_{J^1}(2))$ contains $\mathcal{B}$ and by definition it is exactly $J^1$.  However,
\begin{equation*}
J^1 \cap \PP^4_a = X_a,
\end{equation*}
so if $X_a \subsetneq \mathcal{B}$, then we have a contradiction.  Clearly, if there are only 2
quadrics not containing $\PP^4_a$, then $\mathcal{B}$ has dimension 2 and that is too much.
Finally, if there are 3 quadrics, then $\mathcal{B}$ is a curve of degree 8 containing $X_a$ but
$X_a$ is of degree 6, so there is some extra ``stuff'', a residual conic: a contradiction as well.
\end{proof}

Therefore, the restriction of the dual map has a $\PP^3$-target:
\begin{equation} \label{eq:calD'}
\calD'= \calD|_{\PP^4_a}:\ \PP^4_a \dasharrow \PP^3 = H^0(\PP^4_a,\calI_{X_a}(2))^*.
\end{equation}

So we have reduced the problem to proving that $\Psi = \calD'$, which will be derived from the
following proposition.

\begin{prop} \label{prop:baselocus}
  The linear system of quadrics defining $\Psi$ has base locus $X_a$.
\end{prop}

\begin{proof}
  Certainly, by the way we defined $\Psi$ in Proposition \ref{prop:Psi}, its base locus $B$ is
  contained in $X_a$.  Recall that $X_a$ is embedded into $\PP^4_a = \PP H^0(X_a,\vartheta_a{}^3)^*$
  by
\begin{equation*}
p \mapsto \{ s \in H^0(X_a,\vartheta_a{}^3) \st s(p)=0 \}.
\end{equation*}
So for a point $x \in H^0(\vartheta_a{}^3)^*$, to be in the embedded $X_a$ means that there exists a
(unique) point $p$ on $X_a$ such that
\begin{equation*}
x = ev_p, \text{ i.e. evaluation of sections of $H^0(X_a,\vartheta_a{}^3)$ at $p$.}
\end{equation*}
So let us fix $x \in X_a \subset \PP^4_a$, and its corresponding point $p$ on $X_a$.  We want to see
that $g'(x)$ never has maximal rank.  Remembering the exact diagram \eqref{eq:commdiag3}, elements
of $H^0(X_a,\vartheta_a{}^2 \otimes L^*)$ can be seen as linear combinations of tensors
\begin{equation*}
s \otimes \sigma \in H^0(X_a,\vartheta_a{}^3) \otimes H^0(X_a,\vartheta_a \otimes L)
\end{equation*}
such that $s\cdot\sigma \in H^0(X_a,\vartheta_a{}^4 \otimes L)$ is zero, i.e.
\begin{equation} \label{eq:s.sigma=0}
\forall q \in X_a,\ s(q)\sigma(q)=0.
\end{equation}
But on basic tensors, $g'(x)$ acts as follows
\begin{equation*}
\xymatrix@R=0pt{
  g'(x):\ H^0(X_a,\vartheta_a{}^2 \otimes L^*) \ar[r] & **[r] H^0(X_a,\vartheta_a \otimes L) \\
  s \otimes \sigma \ar@{|->}[r] & **[r] s(p) \cdot \sigma,}
\end{equation*}
and we then see that, because of \eqref{eq:s.sigma=0}, the image of $g'(x)$ is the subspace of
sections in $H^0(X_a,\vartheta_a\otimes~L)$ vanishing at $p$.  But $\vartheta_a\otimes~L$ was
assumed to be globally generated (see proof of Proposition \ref{prop:Psi}), and this subspace is a
proper subspace, i.e. $g'(x)$ degenerates.  Thus we see from \eqref{eq:g'(x) degenerates} that $X_a
\subset B$.
\end{proof}

\begin{prop}  \label{prop:Psi=D'}
  The rational maps $\Psi$ (Proposition \ref{prop:Psi}) and $\calD'$ \eqref{eq:calD'} are equal.
\end{prop}

\begin{proof}
  By Proposition \ref{prop:baselocus}, $\Psi$ is given by a linear subseries of $|\calI_{X_a}(2)|$
  of dimension~3.  But $\calD'$ is given by the complete linear series $|\calI_{X_a}(2)|$ which also
  has dimension~3.  So they are equal.
\end{proof}

Although this last proposition gives a good interpretation of the restricted dual map $\calD'$ from
a vector bundle standpoint, a direct examination of $\calD'$ yields:

\begin{prop} \label{prop:conicfiber-dominant}
  The general fiber of the rational map
\begin{equation*}
\calD':\ \PP^4_a \dasharrow |\calI_{X_a}(2)|^* = \PP^3_a
\end{equation*}
is a conic curve which is a 4-secant of $X_a$.  In particular, $\calD'$ and therefore $\Psi$ are
dominant.
\end{prop}

\begin{proof}
  A point of the target space is a hyperplane (of dimension 3) in $H^0(\PP^4_a,\calI_{X_a}(2))$ (of
  dimension 4), so it is the span of three linearly independent vectors, i.e. 3 quadrics in
  $\PP^4_a$ all containing $X_a$.  In general, their intersection is a degree 8 curve $C$ having
  $X_a$, which is of degree 6, as an irreducible component.  Therefore, the residual curve is a
  conic $Q$.  By the adjunction formula, we see that the $p_a(C)=5$.  Since $g(X_a)=2$ and $g(Q)=0$,
  it follows that $Q$ intersect $X_a$ at 4 points.
\end{proof}

Proposition \ref{prop:Psi=D'} has a rather strong consequence.  One natural question about the
construction of $\Psi$ is: what if we see $x \in \PP^4_a$ as an element of $\PP^4_b$ another
$\PP^4$?  Does the construction give the "same" vector bundle?  From the definition it is not clear
a priori.  But since $\Psi$, or rather $\Psi_a$, coincides with $\calD'$, it does not matter which
$\PP^4$ we choose to define $\Psi(x)$.  So we can globalize the definition and rename $\Psi$ and
$\calD'$.

\begin{cor} \label{cor:global} 
  We have a "global" equality of dominant rational maps:
\begin{equation*}
\Psi = \calD':\ \bigcup_{a\in J} \PP^4_a \dasharrow \Sigma,
\end{equation*}
and
\begin{equation*}
\Sigma \subset \text{Sing}(\check{\calC}_3).
\end{equation*}
\end{cor}

\begin{proof}
  All there is left to prove in the first assertion is that the target space is $\Sigma =
  \text{Sing}(\calC_6)$ and the map is dominant.  From the point of view of $\Psi$ it is clear since
  the construction produced a strictly semi-stable vector bundle. Furthermore, by Proposition
  \ref{prop:conicfiber-dominant}, $\Psi$ is dominant onto the union of $\PP^3_a$.  If we write
  $\vartheta = \OO_X(\Theta|_X)$, we see that
\begin{equation*}
\forall a \in J,\quad \vartheta_a = \vartheta \otimes L_{-a}
\end{equation*}
where $L_{-a}$ is the line bundle on $X$ associated to $-a \in J$.  Therefore, following the
notation of \eqref{eq:Dnu}, we rewrite
\begin{equation*}
\PP^3_a = D \circ \nu(\SU_X(2,\omega_X \otimes \vartheta^{-1} \otimes L_a)),
\end{equation*}
and their union is
\begin{equation*}
\bigcup_{a \in J} \PP^3_a = D \circ \nu(\U_X(2,0)) = \Sigma,
\end{equation*}
and that proves that $\Psi$ is dominant.  From Proposition \ref{prop:conicfiber-dominant}, we also
know that the map $\calD'$ has one-dimensional fibers.  By a property of dual varieties, if some
subvariety of $\calC_3$ not containing the singular locus $J^1(X)$ is mapped onto a lower
dimensional subvariety, then this lower dimensional subvariety lies in the indeterminacy locus of
the inverse dual map, i.e. in the singular locus of the dual variety $\check{\calC}_3$.
\end{proof}

\section{Finishing the proof of the duality and non-abelian Torelli}

\begin{thm} \label{thm:duality}
  The Coble hypersurfaces $\calC_3$ and $\calC_6$ are dual.
\end{thm}

Before we prove the theorem, we state a general auxiliary lemma that we will use in the proof.

\begin{lem} \label{lem:k geq d}
  Let $\mathcal{G}$ be a group acting on the projective space $\PP^n$.  Let $V_d$ be a
  $\mathcal{G}$-invariant hypersurface of degree $d$.  Let $W_k$ be a hypersurface of degree $k$
  such that the intersection $Y=V_d \cap W_k$ is $\mathcal{G}$-invariant.  If $k<d$, then $W_k$ is
  $\mathcal{G}$-invariant.
\end{lem}

\begin{proof}
  Let $F_d$ (resp. $G_k$) be a homogeneous polynomial defining $V_d$ (resp. $W_k$).  Then the
  homogeneous ideal of $Y$ is generated by $F_d$ and $G_k$.  This ideal is $\mathcal{G}$-invariant,
  by this we mean that each homogeneous part of the ideal is an invariant subspace of the vector
  space of homogeneous polynomials of fixed degree.  If $k<d$, then $G_k$ is the only form of degree
  $k$, so it has to be $\mathcal{G}$-invariant.
\end{proof}

To use this lemma, we need a group action.  Let us recall that $\PP^8=|3\Theta|$ and as such there
is an action of $J_3 = (\ZZ/3\ZZ)^4$, the group of $3$-torsion points of the Jacobian of $X$, on
$|3\Theta|$.  This action lifts to an action of the central extension of this group by
$\boldsymbol{\mu_3}$, the group of cubic roots of 1, on the linear space $H^0(J^1,\OO(3\Theta))$.
This central extension called the discrete Heisenberg group.  More details on the Heisenberg group
in the context of the Coble sextic and Coble cubic can be found in \cite{Bea03}, \cite{Ngu},
\cite{NgRa}, \cite{Ort03}.

\begin{proof}[Proof of Theorem \ref{thm:duality}]
  From \cite{Las96}, we know the local structure of the singular locus $\Sigma'$ of $\SU_X(3)$.
  Locally analytically around a general strictly semi-stable bundle $E=F \oplus \det(F)^*$,
  $\SU_X(3)$ looks like a rank-4 quadric in $\A^9$.  Therefore, when we intersect this with a
  general $\A^3$ (through the origin), we get an $A_1$-singularity, i.e a quadratic cone with
  ordinary double point.  Moreover, a standard Chern class computation shows that the the degree of
  the dual variety $\check{\calC}_3$ of $\calC_3$ is a sextic hypersurface (see \cite{Ort03},
  Section 2.4):
\begin{equation*}
\deg(\check{\calC}_3) = 6.
\end{equation*}
Recall first from Corollary \ref{cor:global} that
\begin{equation} \label{eq:sigma in intersection}
\Sigma = \text{Sing}(\calC_6) \subset \text{Sing}(\check{\calC}_3).
\end{equation}
Let us assume that $\calC_6$ and $\check{\calC}_3$ are different.  We write
\begin{equation*}
Y = \calC_6 \cap \check{\calC}_3.
\end{equation*}
Since $\calC_6$ has a singular locus $\Sigma$ of codimension 2, it is irreducible.  Moreover,
$\calC_3$ is an irreducible hypersurface of $\check{\PP}^8$, so its dual variety $\check{\calC}_3$
is also irreducible (see e.g. \cite{GKZ94} Proposition 1.3).  Therefore, $Y$ is connected (by
\cite{FH79} Proposition 1).  We will divide the possible situations into separate cases. \\

\noindent
{\bf Case 1:} $Y$ is reduced, i.e. it has a reduced component.  If we intersect the whole thing with
a general $\PP^3$, we see that the surface $S=\calC_6\cap\PP^3$ has $A_1$-singularities.  Then let's
assume that $T=\check{\calC}_3\cap\PP^3$ is different from $S$.  We denote by $D$ the Cartier
divisor of $S$ defined by the complete intersection with $T$, it is then a connected reduced (not
necessarily irreducible) curve.  In particular, it has no embedded components.  We resolve the 45
rational double points and get 45 exceptional $(-2)$-curves $E_1, \dots, E_{45}$:
\begin{equation*}
\xymatrix{ \coprod_{i=1}^{45} E_i \ar[r] \ar[d] & \tilde{S} \ar[r] \ar[d]^{\pi} & \tilde{\PP}^3 \ar[d]^{\pi} \\
  \{ \text{45 points} \} \ar[r] & S \ar[r] & \PP^3. }
\end{equation*}
We will also denote by $H$ the pullback of the hyperplace section of $S$, therefore its
self-intersection in $\tilde{S}$ is
\begin{equation*}
\int_{\tilde{S}} H^2 = 6.
\end{equation*}
The proper transform $\tilde{D}$ of $D$ under the blowup map $\pi$ is linearly equivalent to
\begin{equation*}
\tilde{D} = 6H - \sum_{i=1}^{45} a_i E_i, \quad a_i \geq 2,
\end{equation*}
since $T$ is also singular at those 45 points, from Corollary \ref{cor:global}.  Recall that $\pi$
is a crepant resolution, so
\begin{align*}
  \omega_{\tilde{S}} &= \pi^*(\omega_{\PP^3} \otimes \OO_{\PP^3}(S) \otimes \OO_S), \\
  &= \OO(2H).
\end{align*}
By the adjunction formula, we can compute the arithmetic genus $p_a(\tilde{D})$ of $\tilde{D}$,
knowing that $\tilde{D}$ is reduced:
\begin{align*}
  2p_a(\tilde{D}) - 2  &= \int_{\tilde{S}} (K_{\tilde{S}} + \tilde{D}) \cdot \tilde{D}, \\
  &= \int_{\tilde{S}} \left(8H - \sum_{i=1}^{45} a_i E_i\right) \cdot \left(6H - \sum_{i=1}^{45} a_i E_i\right), \\
  &= 8 \cdot 6 \cdot 6 + \sum_{i=1}^{45} (-2)a_i^2, \\
  &\leq 288 - 360 = -72,
\end{align*}
because $a_i \geq 2$.  Therefore we see that
\begin{equation*}
p_a(\tilde{D}) \leq -35,
\end{equation*}
which is not possible, since $\tilde{D}$ is an effective Cartier divisor on a nonsigular surface,
therefore without embedded components.  So $S=T$.  Finally, since the intersecting $\PP^3$ was
general, $\calC_6 = \check{\calC}_3$.\\

\noindent
{\bf Case 2:} $Y$ is not reduced.  We write the decomposition of $Y$ into irreduclible component:
\begin{equation*}
Y = a_1 Y_1 + a_2 Y_2 + \dots + a_m Y_m,
\end{equation*}
where $a_i \geq 2$ and $Y_i$ are prime divisors.  But $Y$ is of degree 36, therefore
\begin{equation*}
\sum_{i=1}^m a_i d_i = 36,
\end{equation*}
where $d_i$ denotes the degree of $Y_i$.  Since $\dim \calC_6 \geq 3$, we know by Lefschetz's
Theorem that
\begin{equation*}
\Pic \calC_6 = \ZZ.
\end{equation*}
So each prime divisor $Y_i$ is cut out by a hypersurface in $\PP^8$ and it follows that 6 divides
$d_i$.  Thus the only possible cases are
\begin{itemize}
\item $m=1$: $(a_1,d_1)=(2,18)$ or $(3,12)$ or $(6,6)$.
\item $m=2$: $\{(a_1,d_1),(a_2,d_2)\} = \{(2,12),(2,6)\}$ or $\{(3,6),(3,6)\}$.
\item $m=3$: $\{(a_1,d_1),(a_2,d_2)\,(a_3,d_3)\} = \{(2,6),(2,6),(2,6)\}$.
\end{itemize} 
In every case, we can see that the $a_i$ have a common divisor that is either 2, 3 or 6.  So we can
rewrite
\begin{equation*}
Y = 2Z \text{ or } 3Z \text{ or } 6Z.
\end{equation*}
We know that both $\calC_6$ and $\check{\calC}_3$ are Heisenberg-invariant (see e.g. \cite{Ngu}), so
$Y$ is Heisenberg-invariant, therefore $Z$ is Heisenberg-invariant.  In the cases where $Y=3Z$ or
$6Z$, it means that $Z$ is a quadric or a hyperplane, but there are no Heisenberg-invariant quadrics
or hyperplanes.  But by Lemma \ref{lem:k geq d}, we get a contradiction.  So we are left with one
case: $Y=2Z$ and $Z$ is cut out by a Heisenberg-invariant cubic.  We will also rule this case out.
Again, from \eqref{eq:sigma in intersection}, it follows that $\Sigma \subset Z$.  The involution
$\tau$ of $J^1$ given by
\begin{equation*}
\tau:\ L \mapsto \Omega_X \otimes L^{-1}
\end{equation*}
can be lifted to an involution $\tau$ on $H^0(J^1,\OO(3\Theta))$ (by pulling back sections) and on
$\PP^8=|3\Theta|$.  Then the fixed locus of $\tau$ in $\PP^8$ consists of 2 disjoint projective
spaces a $\PP^4$ and a $\PP^3$.  We know (see \cite{NgRa}, Section 4) that
\begin{equation*}
\Sigma \cap \PP^4 = 2H \cup \{ (15_3)\text{-configuration of lines and points} \},
\end{equation*}
where $H$ is a hyperplane of $\PP^4$.  The configuration of lines and points does not lie in a
hyperplane, so we see that $\Sigma \cap \PP^4$ cannot lie in a cubic hypersurface of $\PP^4$, which
shows that this special $\PP^4$ must lie in the Heisenberg-invariant cubic hypersurface cutting out
$Z$.  But it is easy to check that this special $\PP^4$ cannot be contained in a
Heisenberg-invariant cubic.\footnote{There is a five dimensional space of Heisenberg-invariant
  cubic.  A basis can be found in \cite{Ngu} or \cite{Ort03} and the computations are easy to
  check.}  Contradiction.
\end{proof}

Once this duality is established, we can recover a few well known dualities in classical algebraic
geometry \cite{Ngu} and reinterpret the results in terms of vector bundles.  However, an easy
corollary is the following non-abelian Torelli theorem:

\begin{cor}[Non-abelian Torelli Theorem] \label{cor:torelli}
  Let $X$ and $X'$ be two smooth projective curves of genus 2.  If $\SU_X(3)$ is isomorphic to
  $\SU_{X'}(3)$, then $X$ is isomorphic to $X'$.
\end{cor}

\begin{proof}
  Starting from $\SU_X(3)$, there is a canonical way to retrieve $X$.  We first take the ample
  generator $\Thgen$ of $\Pic(\SU_X(3))$, look at the map associated to the line bundle.  The branch
  locus of the 2-1 map has dual variety a cubic hypersurface in $\PP^8$ singular exactly along the
  the principally polarized Jacobian $(J^1(X),\Theta)$, which determines $X$.
\end{proof}

\begin{rmk}
  The non-abelian Torelli question has been raised ever since the construction of the moduli spaces
  $\SU_X(r,d)$, for any curve $X$ of genus $g\ge 2$.  In \cite{MN68}, D. Mumford and P. E. Newstead
  prove the theorem in the rank $r=2$ case and odd degree determinant ($d=1$), for any $g\ge 2$.  It
  is further generalized to the cases in which $(r,d)=1$ in \cite{NR75} and \cite{Tyu74}.  In all of
  those cases, the moduli spaces are smooth.  Then in \cite{Bal90}, V. Balaji proves the theorem for
  $r=2$, $d=0$ on a curve of genus $g\ge 3$, before A. Kouvidakis and T. Pantev extend the result to
  any $r$ and $d$, still for $g\ge 3$ \cite{KP95}.
  
  In \cite{HR03}, J.-M. Hwang and S. Ramanan introduce a stronger non-abelian Torelli result.  If we
  denote by $\SU_X(r,d)^s$ the moduli space of stable vector bundles, which is open in $\SU_X(r,d)$,
  then $\SU_X(r,d)^s \cong \SU_{X'}(r,d)^s$ implies that $X\cong X'$, and this for any $r$ and $d$,
  but for $g\ge 4$.
  
  Our version of non-abelian Torelli can be shown to be strong, because there still is the ample
  generalized theta divisor on $\SU_X(3)^s$.  Therefore $\SU_X(3)^s$ still dominates $\PP^8 \cong
  |3\Theta|$ and we get the ``open'' Coble sextic $\calC_6-\Sigma$.  The dual map is just the same
  map, as it is defined away from the singular locus.  The closure of the image is the full Coble
  cubic $\calC_3$, therefore we get the desired result, which is a new case for non-abelian Torelli.
\end{rmk}



\end{document}